\documentclass{amsart}
\usepackage{amsmath}
\usepackage{amsfonts,amssymb,amsthm}
\usepackage{mathscinet}
\usepackage{color,comment}
\makeatletter
\newtheorem{theorem}{Theorem}[section]
\newtheorem{lemma}[theorem]{Lemma}
\newtheorem{proposition}[theorem]{Proposition}
\newtheorem{corollary}[theorem]{Corollary}

\theoremstyle{definition}

\newtheorem{definition}[theorem]{Definition}

\newtheorem{problem}[theorem]{Problem}

\newcommand{\ext}{\mathop{\mathrm{ext}}\nolimits}

\newcommand{\Sp}{\mathop{\mathrm{sp}}\nolimits}
\newcommand{\rank}{\mathop{\mathrm{rank}}\nolimits}

\newcommand{\A}{\mathfrak{A}}
\newcommand{\M}{\mathfrak{M}}
\newcommand{\uA}{\widetilde{\mathfrak{A}}}

\newcommand{\cB}{\mathcal{B}}

\newcommand{\cI}{\mathcal{I}}
\newcommand{\cJ}{\mathcal{J}}
\newcommand{\cK}{\mathcal{K}}
\newcommand{\cH}{\mathcal{H}}

\newcommand{\cP}{\mathcal{P}}
\newcommand{\cS}{\mathcal{S}}

\newcommand{\perpn}{\perp_{BJ}}
\newcommand{\perps}{\perp_{BJ}^s}
\newcommand{\perpq}{\perp_{BJ}^q}

\begin{document}
\title[Three types of Birkhoff-James orthogonality]{On relationship among three types of Birkhoff-James orthogonality}

\author[S. Daptari]{Soumitra Daptari}
\address[S. Daptari]{Katsushika Division, Institute of Arts and Sciences, Tokyo University of Science, Tokyo 125-8585, Japan}
\email{daptarisoumitra@rs.tus.ac.jp; daptarisoumitra@gmail.com}

\author[K. Igarashi]{Koki Igarashi}
\address[K. Igarashi]{Department of Applied Mathematics, Graduate School of Science, Tokyo University of Science, Tokyo 162-8601, Japan}
\email{1424501@ed.tus.ac.jp}

\author[J. Nakamura]{Jumpei Nakamura}
\address[J. Nakamura]{Department of Applied Mathematics, Graduate School of Science, Tokyo University of Science, Tokyo 162-8601, Japan}
\email{1425526@ed.tus.ac.jp}

\author[R. Tanaka]{Ryotaro Tanaka}
\address[R. Tanaka]{Katsushika Division, Institute of Arts and Sciences, Tokyo University of Science, Tokyo 125-8585, Japan}
\email{r-tanaka@rs.tus.ac.jp}

\thanks{The fourth named author was supported by JSPS KAKENHI Grant Number JP24K06788. The research of the first named author is sponsored by Tokyo University of Science Post-Doctoral Fellowship under the supervision of the fourth named author.}


\date{}
\begin{abstract}
In this paper, we study three types of Birkhoff-James orthogonality in Hilbert $C^*$-modules, that is, the strong, quasi-strong, and original Birkhoff-James orthogonality. In general, the strong Birkhoff-James orthogonality is stronger than the quasi-strong Birkhoff-James orthogonality, and the quasi-strong Birkhoff-James orthogonality is stronger than the original Birkhoff-James orthogonality. Meanwhile, each reverse implication in this chain requires additional conditions. As the main results, we show that the strong and quasi-strong Birkhoff-James orthogonality are equivalent in a full Hilbert $C^*$-module if and only if the underlying $C^*$-algebra is commutative, and that the equivalence of the quasi-strong and original Birkhoff-James orthogonality in a full Hilbert $C^*$-module implies the primeness of the underlying $C^*$-algebra. Moreover, two examples, explaining the complexity of conditions for full Hilbert $C^*$-modules in which the quasi-strong and original Birkhoff-James orthogonality are equivalent, are given in the $C^*$-algebra settings.
\end{abstract}
\maketitle
\section{Introduction}

This paper is concerned with three types of Birkhoff-James orthogonality in Hilbert $C^*$-modules. The definition of the original Birkhoff-James orthogonality is the following.
\begin{definition}\label{dfBJ}

Let $X$ be a normed space over $\mathbb{K}$, and let $x,y \in X$. Then, $x$ is said to be \emph{Birkhoff-James orthogonal} to $y$, denoted by $x \perpn y$, if $\|x+\lambda y\| \geq \|x\|$ for each $\lambda \in \mathbb{K}$.
\end{definition}
The relation $\perpn$ is recognized as a generalized orthogonality relation since it is equivalent to the usual inner product orthogonality $\perp$ in inner product spaces. The first formulation of this orthogonality dates back to Birkhoff~\cite{Bir35}, and its basics was completed by Day~\cite{Day47} and James~\cite{Jam47a,Jam47b} later, mainly in 1940s. The importance of the Birkhoff-James orthogonality comes from the strong connection with geometric structures of normed spaces. For example, we know characterizations of rotundity and smoothness of norms in terms of the Birkhoff-James orthogonality. Moreover, inner product spaces are characterized as normed spaces in which the relation $\perp_{BJ}$ is symmetric. The readers interested in the Birkhoff-James orthogonality are referred to the surveys by Alonso-Martini-Wu~\cite{AMW12,AMW22} which give good expositions of the existing results on Birkhoff-James orthogonality in normed spaces.

In view of classification, the Birkhoff-James orthogonality gives enough information for identifying normed spaces linearly. To be precise, if a pair of normed spaces admits a linear Birkhoff-James isomorphism, that is, a linear bijection preserving the Birkhoff-James orthogonality in both directions, then they are isometrically isomorphic. This fact immediately follows from better results by Koldobsky~\cite{Kol93} (for the real case) and Blanko-Turn\v{s}ek~\cite{BT06} (for the general case). Based on this situation, recently, nonlinear classification of normed spaces using Birkhoff-James isomorphisms has been begun, for example, in \cite{AGKRZ23,IT22,Tan22a,Tan22b}; see also~\cite{AGKZ22} for understanding the topic. In particular, in the context of $C^*$-algebras, it was shown in \cite{Tan23c} that commutative $C^*$-algebras can be completely classified in $C^*$-algebras by their Birkhoff-James orthogonality structure. Moreover, in \cite{KS25a,KS25b,KSt25}, finite-dimensional $C^*$-algebras were classified in terms of the Birkhoff-James orthogonality, and Birkhoff-James isomorphisms between finite-dimensional $C^*$-algebras were determined. As these examples indicate, the Birkhoff-James orthogonality is getting more important for analyzing $C^*$-algebras as Banach spaces.

The strong-version of Birkhoff-James orthogonality was introduced in Aramba\v{s}i\'{c}-Raji\'{c}~\cite{AR14} by using the module structure of Hilbert $C^*$-modules.
\begin{definition}

Let $\A$ be a $C^*$-algebra, let $\M$ be a Hilbert $\A$-module, and let $x,y \in \M$. Then $x$ is said to be \emph{strongly Birkhoff-James orthogonal} to $y$, denoted by $x \perps y$, if $x \perpn ya$ for each $a \in \A$.
\end{definition}

Considering an approximate unit for $\A$, we can realize that $x \perps y$ implies $x\perpn y$. As expected, the relation $\perps$ in a Hilbert $C^*$-modules captures algebraic information of the underlying $C^*$-algebra. For example, it was shown in \cite{KST18} that every linear strong Birkhoff-James isomorphism, that is, every linear bijection preserving strong Birkhoff-James orthogonality in both directions, between von Neumann algebras is a scalar multiple of a $*$-isomorphism, follows by a unitary multiplication. In other words, von Neumann algebras are distinguished in the sense of $*$-isomorphisms by the combination of the linear and strong Birkhoff-James orthogonality structures. A result on $\cB (\cH)$ similar to above, but for the mutual strong Birkhoff-James orthogonality, was obtained in \cite{AGKRZ22}. There also exists a nonlinear classification theory for strong Birkhoff-James isomorphisms. In~\cite{AGKRTZ24}, Hilbert $C^*$-modules over compact $C^*$-algebras were classified by using strong Birkhoff-James isomorphisms, and in \cite{Tan25}, commutative $C^*$-algebras were characterized and classified in terms of the strong Birkhoff-James orthogonality structure.

As a useful characterization of the Birkhoff-James orthogonality, we know the following result by James~\cite{Jam47b}: Let $X$ be a normed space, and let $x,y \in X$. Then, $x \perpn y$ if and only if $\rho (x)=\|x\|$ and $\rho (y)=0$ for some $\rho \in X^*$ with $\|\rho\|=1$. We also have the Hilbert $C^*$-module version of this result; see~Aramba\v{s}i\'{c}-Raji\'{c}~\cite[Theorem 2.7]{AR12} and also~\cite{BS99} for its prototype. Namely, if $X$ is replaced with a Hilbert $\A$-module $\M$, then $x \perpn y$ if and only if $\rho (\langle x,x \rangle)=\|x\|^2$ and $\rho (\langle x,y \rangle)=0$ for some state $\rho$ of $\A$. Moreover, it follows from \cite[Theorem 2.5]{AR14} that $x \perps y$ if and only if $\rho (\langle x,x\rangle )=\|x\|^2$ and $\rho (|\langle x,y\rangle|^2) =0$ for some state $\rho$ of $\A$, where the state $\rho$ can be chosen pure as was indicated in the proof of \cite[Proposition 3.14]{INRT25b}. These results suggest a natural intermediate of the strong and original Birkhoff-James orthogonality.
\begin{definition}

Let $\A$ be a $C^*$-algebra, let $\M$ be a Hilbert $\A$-module, and let $x,y \in \M$. Then, $x$ is said to be \emph{quasi-strongly Birkhoff-James orthogonal} to $y$, denoted by $x\perpq y$, if $\rho (\langle x,x\rangle)=\|x\|^2$ and $\rho (\langle x,y\rangle )=0$ for some pure state of $\A$.
\end{definition}

It follows that $x \perps y$ implies $x \perpq y$, and $x\perpq y$ implies $x\perpn y$ in any Hilbert $C^*$-module; for $|\rho (\langle x,y \rangle)| \leq \rho (|\langle x,y \rangle|^2)^{1/2}$ whenever $\rho$ is a state. The relation $\perpq$ was first introduced in \cite{INRT25b} for $C^*$-algebras to analyze locally parallel sets and their inclusion relation, and the preceding definition is a natural generalization of \cite[Definition 3.13]{INRT25b}. As interesting features of the quasi-strong Birkhoff-James orthogonality in a $C^*$-algebra $\A$, we know that it is equivalent to the strong Birkhoff-James orthogonality if $\A$ is commutative, and to the original Birkhoff-James orthogonality if $\A$ is $*$-isomorphic to $\cK (\cH)$; see~\cite[Remarks 3.16 and 3.17]{INRT25b}.

The existing results on $\perps$, $\perpq$ and $\perpn$ are sometimes analogous but sometimes completely different. To avoid unintentional duplications of results in future and to develop those theories in Hilbert $C^*$-modules properly, it is important to understand their relationship well. The present paper approaches this theme. In this direction, we only know that Hilbert spaces are the only full Hilbert $C^*$-modules in which $\perps$ and $\perpn$ are equivalent; see~\cite[Theorem 3.7]{AR15}.

The main results of the present paper are Theorems~\ref{SQC} and \ref{th3.7}. We show in Theorem~\ref{SQC} that $\perp\A$ is commutative; and in Theorem~\ref{th3.7} that the underlying $C^*$-algebra $\A$ is prime if $\perpq$ and $\perpn$ are equivalent in a full Hilbert $\A$-module $\M$. It is also shown that the converse to Theorem~\ref{th3.7} is false in general. In fact, as is indicated in Theorem~\ref{Ex-Thm}, there are many primitive $C^*$-algebras in which $\perpq$ and $\perpn$ are not equivalent. Meanwhile, we have a non-simple $C^*$-algebra in which $\perpq$ and $\perpn$ are equivalent. This explains the complexity of conditions for full Hilbert $C^*$-modules in which $\perpq$ and $\perpn$ are equivalent.

\section{Notations and Preliminaries}\label{Sec-np}

Let $X$ be a vector space. A subset $C$ of $X$ is said to be \emph{convex} if $(1-t)x+ty \in C$ whenever $x,y \in C$ and $t \in [0,1]$. A convex subset $D$ of a convex set $C$ is called a \emph{face} of $C$ if $x,y \in C$ and $(1-t)x+ty \in D$ for some $t\in (0,1)$ imply $x,y \in D$. Let $x \in C$. Then, $x$ is called an \emph{extreme point} of $C$ if $\{x\}$ is a face of $C$. The set of all extreme points of $C$ is denoted by $\ext (C)$. If $X$ is a normed space, then the \emph{unit ball} of $X$ means $\{ x \in X : \|x\| \leq \}$. The symbol $X^*$ stands for the (continuous) dual space of $X$. A standard textbook of Banach space theory is \cite{Meg98}.

A $C^*$-algebra is a complex Banach algebra $\A$ with a conjugate-linear involution operation $x \mapsto x^*$ such that $\|x^*x\|=\|x\|^2$ for each $x \in \A$. The algebra $\cB (\cH)$ of all bounded linear operators on a Hilbert space $\cH$ is a typical example of a $C^*$-algebra. The algebra $\cK(\cH)$ of all compact operators on $\cH$ forms a closed two-sided ideal of $\cB (\cH)$ which is itself a $C^*$-algebra. 

Let $\A$ be a $C^*$-algebra, and let $x \in \A$. Then, $x$ is said to be \emph{self-adjoint} if $x^*=x$, and \emph{positive}, denoted by $x \geq 0$, if $x=y^2$ for some self-adjoint element $y \in \A$. We can define a natural partial ordering, so-called the L"{o}wner order, for self-adjoint element $x,y \in \A$ by declaring that $x \leq y$ if $y-x \geq 0$. If $\A$ is unital with unit $1$, then $u \in \A$ is said to be \emph{unital} if $u^*u=uu^*=1$. A spectrum $\Sp (x)$ of $x \in \A$ is defined by
\[
\Sp (x)=\{ \lambda \in \mathbb{C} : \text{$x-\lambda1$ is non-invertible in $\A$}\} .
\]
If $\A$ is non-untial, then $\widetilde{A}$ denotes the unitization of $\A$. The spectrum of an element of a non-unital $C^*$-algebra is considered in its unitization. Let $(e_\lambda)_\lambda$ be a net in $\A$. Then, $(e_\lambda)_\lambda$ is called an \emph{approximate unit} for $\A$ if $\lim_\lambda \|x-xe_\lambda \|=0$ for each $x \in \A$. It is known that every $C^*$-algebra contains an increasing approximate unit consisting of positive elements of the unit ball.

An element $\rho \in \A^*$ is called a \emph{state} of $\A$ if $\|\rho\|=1$ and $\rho (a) \geq 0$ whenever $a \geq 0$. The set of all states of $\A$ is denoted by $\cS(\A)$. A \emph{pure} state is an element of $\ext (\cS(\A))$. Let $\cP (\A) = \ext (\cS (\A))$. Each state $\rho$ of $\A$ induces, by the Gelfand-Naimark-Segal construction, a representation $(\cH_\rho,\pi_\rho,\xi_\rho)$, where $\cH_\rho$ is a Hilbert space, $\pi_\rho :\A \to \cB(\cH_\rho)$ is a $*$-homomorphism, and $\xi_\rho$ is a unit vector such that $\overline{\pi_\rho (\A)\xi_\rho} = \cH_\rho$ and $\rho (x)=\langle \pi_\rho (x)\xi_\rho,\xi_\rho \rangle$ for each $x \in \A$. If $\rho$ is pure, then $\pi_\rho (\A)$ acts irreducibly on $\cH_\rho$, in which case, the representation $(\cH_\rho,\pi_\rho,\xi_\rho)$ is said to be \emph{irreducible}. The Kadison transitivity theorem is an important tool for analyzing irreducible representations. There are many celebrated textbooks for the theory of $C^*$-algebras; see, for example, \cite{Dav96, KR97a, Mur90, Ped18, Sak98}.

Let $\A$ be a $C^*$-algebra. Then, an \emph{inner product $\A$-module $\M$} is a right $\A$-module endowed with an $\A$-valued inner product $\langle \cdot, \cdot\rangle:\M \times \M\longrightarrow \A$ satisfying the following conditions:
\begin{enumerate}
    \item $\langle x,\alpha y+\beta z \rangle=\alpha \langle x,y\rangle + \beta \langle x,z\rangle$ for each $x,y,z\in \M$ and each $\alpha,\beta \in \mathbb{C}$.
    \item $\langle x,ya\rangle=\langle x,y\rangle a$ for each $x,y\in \M$ and each $a\in \A$.
    \item $\langle x,y\rangle^*=\langle y,x\rangle$ for each $x,y\in \M$.
    \item $\langle x,x\rangle\geq 0$ for each $x\in \M$, and $\langle x,x\rangle= 0$ implies $x=0$.
\end{enumerate}
The norm of an element $x \in \M$ is defined by $\|x\|=\|\langle x,x \rangle\|^{1/2}$. If $\M$ is complete with respect to this norm, then $\M$ is called a \emph{Hilbert $\A$-module}. A Hilbert $\mathbb{C}$-module corresponds to a Hilbert space. Since $\langle xa,y \rangle = \langle y,xa\rangle^*= a^*\langle x,y\rangle$ and
\[
0 \leq\langle xa,xa \rangle =a^*\langle x,x \rangle a \leq \|x\|^2a^*a,
\]
it follows that $\|xa\|=\|\langle xa,xa \rangle\|^{1/2}\leq \|x\|\|a\|$ for each $x \in \M$ and each $a \in \A$. A Hilbert $\A$-module $\M$ is said to be \emph{full} if the linear span of $\{ \langle x,y \rangle : x,y \in \M \}$ is norm-dense in $\A$. If $\M$ is a full $\A$-module $\M$ and $a \in \A \setminus \{0\}$, then $xa \neq 0$ for some $x \in \M$. We note that $\A$ is itself a full Hilbert $\A$-module with the inner product defined by $\langle x,y\rangle =x^*y$ for each $x,y \in \A$. For the basics of Hilbert $C^*$-modules, the readers are referred to \cite[Chapter 1]{Lan95}.

\section{Main results}\label{Sec-main}

First, we consider the equivalence of $\perps$ and $\perpq$ in full Hilbert $C^*$-modules. The following theorem improves \cite[Proposition 3.3]{AR15}.
\begin{theorem}\label{SQC}

Let $\A$ be a $C^*$-algebra, let $\M$ be a full Hilbert $\A$-module. Then $\perps$ and $\perpq$ are equivalent on $\M$ if and only if $\A$ is commutative.
\begin{proof}
Suppose that $\A$ is commutative.  Let $x,y \in \M$ be such that $x \perpq y$, and let $\rho \in \cP (\A)$ be such that $\rho (\langle x,x \rangle )=\|x\|^2$ and $\rho (\langle x,y\rangle )=0$. Since every pure state of $\A$ is multiplicative, it follows that $\rho (\langle x,ya\rangle )=\rho (\langle x,y\rangle a) = \rho (\langle x,y\rangle)\rho (a)=0$ for each $a \in \A$. Hence, $x \perpn ya$ for each $a \in \A$, that is, $x \perps y$.

Conversely, suppose that $\A$ is non-commutative. We claim that there exists a pair $(x,\rho) \in \M \times \cP(\A)$ such that $\dim \cH_\rho \geq 2$ and $\rho (\langle x,x \rangle )=\|x\|^2=1$. To show this, let $\cP_0(\A)$ be the set of all multiplicative states of $\A$. If $\cP_0(\A) = \emptyset$, then the claim is true. In the case $\cP_0 (\A) \neq \emptyset$, we have a $*$-homomorphism $a \mapsto (\rho (a))_{\rho \in \cP_0(\A)}$ from $\A$ into $\ell_\infty (\cP_0(\A))$. Since $\A$ is non-commutative, there exists a nonzero $a \in \A$ such that $\rho (a)=0$ for each $\rho \in \cP_0(\A)$. By the fullness of $\M$, we obtain $xa \neq 0$ for some $x \in \M$. It may be assumed that $\|xa\|=1$. Let $\tau \in \cP(\A)$ be such that $\tau (\langle xa,xa\rangle )=\|xa\|^2=1$. We note that $\dim \cH_\tau \geq 2$. Indeed, if $\dim \cH_\tau =1$, then $\pi_\tau(a) = \tau (a)I \in \cB (\cH_\rho) \cong \mathbb{C}$, which ensures $\tau \in \cP_0(\A)$. However, this is impossible since $\tau \in \cP_0(\A)$ implies $\tau (a)=0$ and $\tau (\langle xa,xa \rangle )= \tau (\langle xa,x\rangle )\tau (a)$. Therefore, $\dim \cH_\rho \geq 2$. This proves our claim.

Let $(x,\rho) \in \M \times \cP(\A)$ be a pair such that $\dim \cH_\rho \geq 2$ and $\rho (\langle x,x \rangle )=\|x\|^2=1$. Since $\langle \pi_\rho (\langle x,x \rangle )\xi_\rho,\xi_\rho \rangle =\rho (\langle x,x \rangle ) =1$, it follows that $\pi_\rho (\langle x,x \rangle)\xi_\rho = \xi_\rho$. We divide our argument into three cases:

(I) $\rank \pi_\rho (\langle x,x \rangle )=1$. In this case, $\pi_\rho (\langle x,x\rangle )=E_{\xi_\rho}$. Let $\xi$ be a unit vector in $\cH_\rho$ such that $\xi_\rho \perp \xi$. By the Kadison transitivity theorem, there exist a self-adjoint element $a \in \A$ and a unitary element $u \in \uA$ such that $\pi_\rho (a)\xi_\rho = \xi_\rho$, $\pi_\rho (a)\xi=\xi$, and $\pi_\rho (u)\xi_\rho=\xi$. Replacing $h$ with $|\sin (2^{-1}\pi a)|$ if necessary, we may assume that $a$ is positive and $\|a\|=1$. It follows that $\|xa\|^2=\|\langle xa,xa\rangle\|=\|a\langle x,x\rangle a\| \leq 1$,
\[
\rho (\langle xa,xa \rangle ) = \rho (a\langle x,x \rangle a)=\langle \pi_\rho (\langle x,x \rangle )\pi_\rho (a)\xi_\rho,\pi_\rho (a)\xi_\rho \rangle =1
\]
and
\begin{align*}
\rho (\langle xa,xau \rangle ) = \rho (a \langle x,x \rangle au ) 
&= \langle \pi_\rho (\langle x,x \rangle )\pi_\rho (a)\pi_\rho (u)\xi_\rho ,\pi_\rho (a)\xi_\rho \rangle \\
&= \langle E_{\xi_\rho}\xi, \xi_\rho \rangle \\
&= 0 .
\end{align*}

(II) $\rank \pi_\rho (\langle x,x \rangle) \geq 2$ and $\Sp (\pi_\rho (\langle x,x\rangle )) \cap (0,1) =\emptyset$. In this case, $E=\pi_\rho (\langle x,x \rangle )$ is a projection with $\xi_\rho \in E(\cH_\rho)$. Let $\xi$ be a unit vector in $E(\cH_\rho)$ such that $\xi_\rho \perp \xi$. The same pair $(a,u)$ as (I) satisfies $\|xa\|^2=\rho (\langle xa,xa \rangle)=1$ and
\[
\rho (\langle xa,xau\rangle ) = \langle \pi_\rho (\langle x,x \rangle)\pi_\rho (a)\pi_\rho (u)\xi_\rho ,\pi_\rho(a)\xi_\rho \rangle = \langle \xi,\xi_\rho \rangle =0 .
\]

(III) $\rank \pi_\rho (\langle x,x \rangle) \geq 2$ and $\Sp (\pi_\rho (\langle x,x\rangle )) \cap (0,1) \neq \emptyset$. By functional calculus, we obtain a pair of positive norm-one elements $h,k \in \A$ commuting with $\langle x,x \rangle$ such that $hk=0$, $\pi_\rho (h)\xi_\rho =\xi_\rho$ and $\pi_\rho (\langle x,x \rangle  k)\neq 0$. Let $\xi \in \cH_\rho$ be such that $\pi_\rho (\langle x,x \rangle )\pi_\rho (k)\xi\neq 0$. It may be assumed that $\|\pi_\rho (k)\xi\|=1$. We note that $\pi_\rho (\langle x,x \rangle )\pi_\rho (k)\xi = \pi_\rho (k)\pi_\rho (\langle x,x \rangle )\xi$ and
\[
\langle \xi_\rho ,\pi_\rho (\langle x,x \rangle)\pi_\rho (k)\xi\rangle = \langle \pi_\rho (h)\xi_\rho,\pi_\rho (k)\pi_\rho (\langle x,x\rangle)\xi \rangle = 0 = \langle \xi_\rho,\pi_\rho (k)\xi\rangle
\]
By the Kadison transitivity theorem, as in (I), we can find a positive element $a \in \A$ and a unitary element $u \in \uA$ such that $\|a\|=1$, $\pi_\rho (a)\xi_\rho = \xi_\rho$, $\pi_\rho (a)\pi_\rho (k)\xi=\pi_\rho (k)\xi$, and $\pi_\rho (u)\xi_\rho = \pi_\rho (k)\xi$. It follows that $\|xa\| \leq 1$, 
\[
\rho (\langle xa,xa \rangle ) = \rho (a\langle x,x \rangle a)=\langle \pi_\rho (\langle x,x \rangle )\pi_\rho (a)\xi_\rho,\pi_\rho (a)\xi_\rho \rangle =1
\]
and
\begin{align*}
\rho (\langle xa,xau \rangle ) = \rho (a \langle x,x \rangle au ) 
&= \langle \pi_\rho (\langle x,x \rangle )\pi_\rho (a)\pi_\rho (u)\xi_\rho ,\pi_\rho (a)\xi_\rho \rangle \\
&= \langle \pi_\rho (\langle x,x \rangle ) \pi_\rho (k)\xi,\xi_\rho \rangle \\
&= 0 .    
\end{align*}
Thus, in either case, we have a norm-one element $x \in \M$, a positive norm-one element $a \in \A$ and a unitary element $u \in \uA$ such that $xa \perpq xau$. Meanwhile, we derive $xa \not \perps xau$ by $0 \leq a \leq 1$ and
\[
\|xa-(xau)u^*a\| = \|x(a-a^2)\| \leq \|a-a^2\| = 1/4<\|xa\| .
\]
This completes the proof.
\end{proof}
\end{theorem}
Next, we consider the equivalence of $\perpq$ and $\perpn$ in full Hilbert $C^*$-modules. In this direction, we have a very strong sufficient condition and an algebraic necessary conditon for the equivalence. The following proposition gives a former one.
\begin{proposition}\label{algebra-to-module}

Let $\A$ be a $C^*$-algebra, and let $\M$ be a Hilbert $\A$-module. If $\perpq$ and $\perpn$ are equivalent on $\A$, then so is on $\M$. In particular, $\perpq$ and $\perpn$ are equivalent on every Hilbert $\cK(\cH)$-module.
\begin{proof}
Suppose that $a \perpq b$ if and only if $a \perpn b$ for each $a,b \in \A$. It follows that $\langle x,x \rangle \perpq \langle x,y \rangle$ if and only if $\langle x,x\rangle \perpn \langle x,y \rangle$ for each $x,y \in \M$. We also know that $x \perpn y$ is equivalent to $\langle x,x \rangle \perpn \langle x,y \rangle$ for each $x,y \in \M$. Thus, $\perpq$ and $\perpn$ are equivalent on $\M$. The last statement follows from \cite[Remark 3.17]{INRT25b}.
\end{proof}
\end{proposition}
Next, we consider a necessary condition. Recall that a $C^*$-algebra $\A$ is said to be \emph{prime} if $\cI \cap \cJ \neq \{0\}$ for each pair of nontrivial closed two-sided ideals $\cI$ and $\cJ$ of $\A$. The following lemma is essentially a version of \cite[Lemma I.9.15]{Dav96}.
\begin{lemma}\label{pure-ideal}

Let $\A$ be a $C^*$-algebra, let $\cI$ be a closed two-sided ideal of $\A$, and let $\rho$ be a pure state of $\A$ such that $\|\rho |\cI\|=1$. Then, $\rho |\cI$ is a pure state of $\cI$.
\begin{proof}
Since $\|\rho \|=\|\rho |\cI\| =1$, by \cite[Proposition II.6.4.16]{Bla06}, $\rho$ is a unique extension of $\rho |\cI$ to $\A$ and $\rho (x) = \lim_a \rho (xe_a)$ for each $x \in \A$, where $(e_a )_a$ is an (increasing) approximate unit. Now, let $\tau_1,\tau_2$ be positive linear functionals on $\cI$ such that $\|\tau_1\|\leq 1$, $\|\tau_2\|\leq 1$ and $(1-t)\tau_1+t\tau_2 = \rho | \cI$. In this case, again by \cite[Proposition II.6.4.16]{Bla06}, $\tau_j$ extends uniquely to a positive linear functional $\rho_j$ on $\A$ satisfying $\|\rho_j\|=\|\tau_j\|$, in which case, $\rho_j (x) = \lim_a \tau_j (xe_a)$. It follows from
\begin{align*}
((1-t)\rho_1+t\rho_2)(x)
&= (1-t)\lim_a \tau_1 (xe_a) + t\lim_a \tau_2 (xe_a) \\
&= \lim_a ((1-t)\tau_1+t\tau_2)(xe_a) \\
&= \lim_a \rho (xe_a) \\
&= \rho (x)
\end{align*}
for each $x \in \A$ that $(1-t)\rho_1+t\rho_2 =\rho$. Since $\rho$ is pure, we obtain $\rho = \rho_1= \rho_2$ and $\rho |\cI = \tau_1 =\tau_2$. This shows that $\rho |\cI$ is a pure state of $\cI$.
\end{proof}
\end{lemma}
\begin{lemma}\label{either}

Let $\A$ be a $C^*$-algebra, let $\cI$ and $\cJ$ be nontrivial closed two-sided ideals of $\A$ such that $\cI \cap \cJ=\{0\}$, and let $\rho$ be a pure state of $\A$ such that $\|\rho |(\cI+\cJ)\|=1$. Then, $\rho |\cI = 0$ or $\rho |\cJ = 0$.
\begin{proof}
First, we note that $\cI+\cJ$ is also a closed two-sided ideal of $\A$ by \cite[Corollary II.5.1.3]{Bla06}. From this, Lemma~\ref{pure-ideal} ensures that $\rho |(\cI+\cJ)$ is a pure sate of $\cI +\cJ$. Let $(e_a)_a$ and $(f_b)_b$ be increasing approximate units for $\cI$ and $\cJ$, respectively. Using~\cite[Technique 2.1.32]{Meg98}, we can find subnets $(e_c)_c$ of $(e_a)_a$ and $(f_c)_c$ of $(f_b)_b$ having the same index set. Since $\cI \cap \cJ=\{0\}$ and $(e_c)_c$ and $(f_c)_c$ are also increasing approximate units for $\cI$ and $\cJ$, it follows that $\|e_c+f_c\| \leq 1$ for each $c$ and
\[
\lim_c (e_c+f_c)(x+y)= \lim_c (e_cx+f_cy) = x+y
\]
for each $x \in \cI$ and each $y \in \cJ$. Hence, $(e_c+f_c)_c$ is an increasing approximate unit for $\cI + \cJ$. Combining this with the fact that $\rho |\cI \geq 0$, $\rho |\cJ \geq 0$ and $\rho |(\cI +\cJ) \geq 0$, we obtain
\[
1=\|\rho|(\cI+\cJ)\| = \lim_c \rho (e_c+f_c) = \lim_c \rho (e_c)+\lim_c \rho (f_c) = \|\rho|\cI\|+\|\rho |\cJ\|
\]
by~\cite[Lemmas I.9.5 and I.9.6]{Dav96}.

Now, let $\tau_1 (x+y) = \rho(x)$ and $\tau_2 (x+y)=\rho(y)$ for each $x+y \in \cI +\cJ$. Then, $\tau_j$ is well-defined by $\cI \cap \cJ=\{0\}$ and linear. If $x+y \in \cI+\cJ$ and $x+y \geq 0$, then $x^*=x$ and $y^*=y$, which implies that $x^3=x(x+y)x \geq 0$ and $y^3 =y(x+y)y \geq 0$. Therefore, $x \geq 0$ and $y \geq 0$ by the spectral mapping theorem. This shows that $\tau_1 (x+y) = \rho (x) \geq 0$ and $\tau_2 (x+y) = \rho (y) \geq 0$. Namely, $\tau_1 \geq 0$ and $\tau_2 \geq 0$. Combining this with~\cite[Lemmas I.9.5 and I.9.6]{Dav96}, we obtain $\|\tau_1\|=\lim_c \tau_1 (e_c+f_c) = \|\rho |\cI\|$ and $\|\tau_2\| = \lim_c \tau_2(e_c+f_c) = \|\rho |\cJ\|$. Finally, if $\|\tau_1\|\|\tau_2\|\neq 0$, then
\[
\rho |(\cI +\cJ) = \tau_1+\tau_2 = \|\tau_1\|(\|\tau_1\|^{-1}\tau_1)+\|\tau_2\|(\|\tau_2\|^{-1}\tau_2) ,
\]
while $\|\tau_1\|^{-1}\tau_1 \neq \|\tau_2\|^{-1}\tau_2$ by $\lim_c (\|\tau_1\|^{-1}\tau_1)(e_c) = 1 \neq 0 =\lim_c (\|\tau_2\|^{-1}\tau_2)(e_c)$. This contradicts the fact that $\rho |(\cI +\cJ)$ is pure. Thus, $\|\tau_1\|\|\tau_2\|=0$, that is, $\rho |\cI =0$ or $\rho |\cJ=0$.
\end{proof}
\end{lemma}
The following result gives an algebraic necessary condition for full Hilbert $C^*$-modules on which $\perpq$ and $\perpn$ are equivalent.
\begin{theorem}\label{th3.7}

Let $\A$ be a $C^*$-algebra, and let $\M$ be a full $\A$-module. If $\perpq$ and $\perpn$ are equivalent on $\A$, then $\A$ is prime.
\begin{proof}
Suppose that $\A$ is not prime. Then, there exists a pair of nontrivial closed two-sided ideals $\cI$ and $\cJ$ of $\A$ such that $\cI \cap \cJ=\{0\}$. Let $a \in \cI$ and $b \in \cJ$ be nonzero. Since $\M$ is full, we have $xa \neq 0$ and $yb \neq 0$ for some $x,y \in \M$. It may be assumed that $\|xa\|=\|yb\|=1$. We note that
\begin{align*}
\langle \alpha xa+\beta yb,\alpha xa+\beta yb\rangle 
&= |\alpha|^2\langle xa,xa \rangle + \alpha \overline{\beta}\langle xa,yb\rangle +\overline{\alpha}\beta\langle yb,xa\rangle +|\beta|^2\langle yb,yb\rangle \\
&= |\alpha|^2\langle xa,xa \rangle+|\beta|^2\langle yb,yb\rangle
\end{align*}
for each $\alpha ,\beta \in \mathbb{C}$ by $a^*\langle x,y \rangle b,b^*\langle y,x \rangle a \in \cI \cap \cJ=\{0\}$. Combining this with 
\[
\langle xa,xa\rangle \langle yb,yb\rangle =a^*\langle x,xa\rangle \langle yb,y\rangle b \in \cI \cap \cJ=\{0\},
\]
we obtain $\|\alpha xa+\beta yb\| = \max \{ |\alpha |,|\beta |\}$ for each $\alpha ,\beta \in \mathbb{C}$. Hence, it follows from
\begin{align*}
\|(xa+yb)+\lambda (xa-yb)\| = \|(1+\lambda)xa+(1-\lambda)yb\| = \max \{ |1+\lambda|,|1-\lambda|\} \geq 1
\end{align*}
for each $\lambda \in \mathbb{C}$ that $xa+yb \perpn xa-yb$. Meanwhile, if $\rho$ is a pure state of $\A$ such that
\[
1=\rho (\langle xa+yb,xa+yb\rangle) = \rho (\langle xa,xa \rangle )+\rho (\langle yb,yb \rangle),
\]
then $\rho |\cI=0$ or $\rho |\cJ=0$ by Lemma~\ref{either}. In the case $\rho |\cJ = 0$, we obtain $\rho (\langle xa,xa \rangle) =1$ and
\begin{align*}
\rho (\langle xa+yb,xa-yb\rangle ) = \rho (\langle xa,xa \rangle ) =1 ;
\end{align*}
while, in the caes $\rho |\cI=0$, we get $\rho (\langle yb,yb\rangle )=1$ and
\begin{align*}
\rho (\langle xa+yb,xa-yb\rangle ) = -\rho (\langle yb,yb \rangle ) =-1 ;
\end{align*}
Therefore, $xa+yb \not \perpq xa-yb$.
\end{proof}
\end{theorem}

\begin{corollary}[Aramba\v{s}i\'{c} and Raji\'{c}, 2015]

Let $\A$ be a $C^*$-algebra, and let $\M$ be a full Hilbert $\A$-module. Then, $\perps$ and $\perpn$ are equivalent on $\M$ if and only if $\A = \mathbb{C}$.
\begin{proof}
If $\perps$ and $\perpn$ are equivalent on $\M$, then Theorems~\ref{SQC} and \ref{th3.7} ensures that $\A$ is commutative and prime. Since every prime $C^*$-algebra has the trivial center, it follows that $\A = \mathbb{C}$. The converse follows from Theorem~\ref{SQC} and Proposition~\ref{algebra-to-module}.
\end{proof}
\end{corollary}


\section{Remarks}\label{Sec-rem}

We conclude this paper with some remarks on Theorem~\ref{th3.7}. It will turn out that the converse to Theorem~\ref{th3.7} is not true and conditions for Hilbert $C^*$-modules on which $\perpq$ and $\perpn$ are equivalent is rather complicated even in the $C^*$-algebra settings.

Before starting, we recall about a dichotomy for pure states. Let $\A$ be a $C^*$-algebra acting on a Hilbert space $\cH$. Suppose that $\cK (\cH) \subset \A$. Then, the following hold:
\begin{itemize}
\item[{\rm (i)}] If $\rho \in \cP (\A)$ is ultraweakly continuous on $\uA$, then $\rho$ is a vector state.
\item[{\rm (ii)}] If $\rho$ is singular on $\uA$, then $\cK (\cH) \subset \ker \rho$.
\item[{\rm (iii)}] If $\rho \in \cP (\A)$, then $\rho$ is a vector state or $\cK (\cH) \subset \ker \rho$.
\end{itemize}
We explain these clauses for the sake of completeness.

(i) We note that $\rho$ extends uniquely to a pure state of $\cP (\uA)$. Moreover, the ultraweak continuity of $\rho$ on $\uA$ ensures that it extends (uniquely) to an ultraweakly continuous state of $\cB (\cH)$ by \cite[Corollary 10.1.11]{KR97b}. From this and \cite[Theorem 7.1.12]{KR97b}, we have an at most countable orthogonal system of nonzero vectors $(\zeta_n)_n \subset \cH$ such that $\sum_{n \in N}\|\zeta_n\|^2=1$ and $\rho = \sum_{n \in N} \omega_{\zeta_n}=\sum_{n \in N}\|\zeta_n\|^2 \omega_{\|\zeta_n\|^{-1}\zeta_n}$ on $\cB (\cH)$. Combining this with $\rho \in \cP (\A)$, we conclude that $N$ is a singleton. Namely, $\rho = \omega_\zeta$ for some unit vector $\zeta \in \cH$.

(ii) We assume that $\rho$ is singular on $\uA$. By \cite[Corollary 10.1.16]{KR97b}, whenever $B,C \in \uA$, the functional $\rho_0$ defined by $\rho_0(A)=\rho (BAC)$ for each $A \in \uA$ is singular or zero. Let $\xi$ be a unit vector in $\cH$. It follows from
\[
\rho (E_\xi  AE_\xi ) = \omega_\xi (A)\rho (E_\xi)
\]
for each $A \in \uA$ that $\rho (E_\xi)\omega_\xi$ is singular as well as ultraweakly continuous. Hence, $\rho (E_\xi)=0$. Thus, $\rho$ vanishes on every rank-one projection on $\cH$, which ensures $\cK (\cH ) \subset \ker \rho$.

(iii) By \cite[Theorem 10.1.15]{KR97b}, $\rho$ (which can be considered as a pure state of $\uA$) is decomposed into the sum of ultraweakly continuous part $\rho_u$ and singular part $\rho_s$, in which case, $\|\rho_u\|+\|\rho_s\|=1$. If $\rho \neq \rho_u$ and $\rho \neq \rho_s$, then $\|\rho_u\|>0$ and $\|\rho_s\|>0$, which implies that
\[
\rho = \|\rho_u\|(\|\rho_u\|^{-1}\rho_u) + \|\rho_s\|(\|\rho_s\|^{-1}\rho_s) .
\]
This contradicts $\rho \in \cP (\A)$. Hence, $\rho = \rho_u$ or $\rho = \rho_s$. The conclusion follows from (i) and (ii).

Now, we present two contrasting results explaining the complexity of conditions for conditions for $C^*$-algebras on which $\perpq$ and $\perpn$ are equivalent.
\begin{theorem}\label{Ex-Thm}

Let $\cH$ be a separable infinite-dimensional Hilbert space, and let $\A$ be a $C^*$-algebra acting on $\cH$. Suppose that $\cK (\cH) \subset \A$. Then, $\perpq$ and $\perpn$ are equivalent on $\A$ if and only if $\A =\cK (\cH)$.
\begin{proof}
It suffices to prove the `only if' part. Suppose to the contrary that $\A \neq \cK (\cH)$. Then, there exists a self-adjoint element $A \in \A$ such that $\|A+\cK(\cH)\|=1$. By \cite[Proposition II.5.1.5]{Bla06}, it may be assumed that $\|A\|=1$. In this case, we also have $\|A^2\|=1=\|A^2+\cK (\cH)\|$. Set $B=A^2 \geq 0$ for short. If $B$ is not norm attaining, then we replace $B$ with $(I-E)B(I-E)+E$ for some rank-one projection $E$. We note that $0 \leq (I-E)B(I-E)+E \leq I$ and
\[
(I-E)B(I-E)+E = B-BE-EB+EBE+E \in B+\cK (\cH) \subset\A .
\]
In what follows, we assume that $BE=E$. Now, let $\xi_1$ be the unit vector such that $E=E_{\xi_1}$, and let $(\xi_n)_n$ be an orthonormal basis for $\cH$. Set
\[
C =B^{1/2}\left(I-\sum_{n \geq 2}\frac{1}{n}E_{\xi_n} \right) B^{1/2} = B-\sum_{n \geq 2}\frac{1}{n}B^{1/2}E_{\xi_n}B^{1/2} \in B+\cK (\cH) 
\]
and $D=C-2E \in \A$. Since $0 \leq C \leq B \leq I$ and $1=\|B+\cK (\cH)\| \leq \|C\|$, it follows that $\|C\|=1$.

We claim that $C \perpn D$. Since $\|C+\cK (\cH)\|=1$, there exists a state $\tau_0$ of $\A /\cK (\cH)$ such that $\tau_0 (C+\cK (\cH)) = 1$. Set $\tau = \tau_0 \circ \pi$, where $\pi:\A \to \A /\cK (\cH)$ is the quotient map. Then, $\tau \in \cS (\A)$, $\tau (C)=1$ and $\cK (\cH) \subset \ker \tau$. It follows from $B\xi_1=B^{1/2}\xi_1=\xi_1$ that
\[
\left( \frac{1}{2}(\omega_{\xi_1}+\tau) \right)(C) = 1
\]
and
\[
\left( \frac{1}{2}(\omega_{\xi_1}+\tau) \right)(D) = 1-\left( \frac{1}{2}(\omega_{\xi_1}+\tau) \right)(2E) = 0 .
\]
Hence, $C \perpn D$.

Next, let $\rho \in \cP (\A)$ be such that $\rho (C)=1$. If $\rho$ is singular, then $\rho (D) = \rho (C) =1$. If $\rho = \omega_\xi$ for some $\xi \in \cH$ with $\|\xi\|=1$, then $\sum_n E_{\xi_n} =I$ implies
\begin{align*}
1 = \langle C\xi,\xi \rangle
&= \left\langle \left( I-\sum_{n\geq 2}\frac{1}{n}E_{\xi_n}\right)B^{1/2}\xi,B^{1/2}\xi \right\rangle \\
&= \langle E_{\xi_1}B^{1/2}\xi,B^{1/2}\xi\rangle+\sum_{n \geq 2}\left( 1-\frac{1}{n}\right)\langle E_{\xi_n}B^{1/2}\xi,B^{1/2}\xi\rangle \\
&= |\langle B^{1/2}\xi,\xi_1\rangle|^2 + \sum_{n \geq 2}\left(1-\frac{1}{n}\right)|\langle B^{1/2}\xi,\xi_n\rangle|^2 \\
&\leq \sum_n |\langle B^{1/2}\xi,\xi_n \rangle|^2 \\
&=\|B^{1/2}\xi\|^2.
\end{align*}
Therefore, $\|B^{1/2}\xi\|=1$, $\langle B^{1/2}\xi,\xi_n\rangle = 0$ whenever $n \geq 2$, and $B^{1/2}\xi=\lambda \xi_1$ for some $\lambda \in \mathbb{C}$ with $|\lambda |=1$. In particular, we obtain $C\xi = \lambda B^{1/2}\xi_1 = \lambda \xi_1$. Meanwhile, it follows from $1=\rho (C) \leq \rho (C^2)^{1/2}=1$ that $\rho ((I-C)^2) = \rho (I-2C+C^2)=0$ and
\[
|\rho ((I-C)T)| = |\rho (T^*(I-C))| \leq \rho (T^*T)^{1/2}\rho ((I-C)^2)^{1/2} = 0
\]
for each $T \in \A$. Thus, $\rho (T) = \rho (CTC) = \omega_\xi (CTC) = \omega_{\xi_1}(T)$ for each $T \in \A$, which together with $E=E_{\xi_1}$ shows that $\rho (D)=\rho (C-2E) =-1 \neq 0$. Consequently, $C \not \perpq D$.
\end{proof}
\end{theorem}
\begin{theorem}\label{th5.3}

Let $\cH$ be a nonseparable infinite-dimensional Hilbert space. Then, $\perpq$ and $\perpn$ are equivalent on $\cK (\cH)^\sim$.
\begin{proof}
We infer that $\cS (\cK (\cH)^\sim)$ contains a unique singular element $\rho_\infty$ which is defined by $\rho_\infty (A+cI)=c$ for each $A+cI \in \cK (\cH)^\sim$. We have $\rho_\infty \in \cP (\cK (\cH)^\sim)$. Let $A,B \in \cK (\cH)$ and $c,d \in \mathbb{C}$ be such that $A+cI\perpn B+dI$, and let $\rho \in \cS (\cK (\cH))$ be such that $\rho (|A+cI|)=\|A+cI\|$ and $\rho ((B+dI)^*(A+cI))=0$. We may assume that $\|A+cI\|=1$. By \cite[Theorem 10.1.15]{KR97b}, we have a pair of ultraweakly continuous and singular functionals $(\rho_u,\rho_s)$ satisfying $\|\rho_u\|+\|\rho_s\|=1$ and $\rho = \rho_u + \rho_s$. Since $1=\rho (I)=\rho_u (I)+\rho_s (I)$, it follows that $\|\rho_u\|=\rho_u (I)$ and $\|\rho_s\|=\rho_s (I)$, which implies that $\rho_u \geq 0$ and $\rho_s \geq 0$. Moreover, the identity
\[
\rho_u (|A+cI|) +\rho_s (|A+cI|) = \|A+cI\|=\|\rho_u\|\|A+cI\|+\|\rho_s\|\|A+cI\|
\]
ensures $\rho_u (|A+cI|) =\|\rho_u\|\|A+cI\|$ and $\rho_s (|A+cI|) = \|\rho_s\|\|A+cI\|$. We note that $\rho_s = \|\rho_s\|\rho_\infty$. Moreover, by \cite[Corollary 10.1.11]{KR97b}, $\rho_u$ extends (uniquely) to an ultraweakly continuous state of $\cB (\cH)$. From this and \cite[Theorem 7.1.12]{KR97b}, we have an at most countable orthogonal system of nonzero vectors $(\zeta_n)_{n \in N} \subset \cH$ such that $\rho_u = \sum_{n \in N} \omega_{\zeta_n} = \sum_{n \in N}\|\zeta_n\|^2\omega_{\|\zeta_n\|^{-1}\zeta_n}$ on $\cB (\cH)$ and $\sum_n \|\zeta_n\|^2=\|\rho_u\|$. In particular, by $\rho_u (|A+cI|) =\|\rho_u\|\|A+cI\|$, we obtain $|A+cI|\zeta_n = \zeta_n$ for each $n \in N$, which means that $|A+cI|E=E$ for the projection $E$ from $\cH$ onto $[\{ \zeta_n : n \in N\}]$.

Suppose that $|c|<1$. In this case, $\rho = \rho_u$ by
\[
\rho_\infty (|A+cI|^2) = |c|^2<1=\|A+cI\|^2 .
\]
We claim that $E$ is finite-rank (or equivalently, $N$ is finite). To show this, suppose to the contrary that $E(\cH)$ contains an orthonormal sequence $(\xi_n)_n$. Since $(\xi_n)_n$ is weakly null, we have $\lim_n \|C\xi_n\|=0$ for each $C \in \cK(\cH)$ and
\[
1 = \|\xi_n\|=\||A+cI|^2\xi_n\| \leq \|(A^*A+cA^*+\overline{c}A)\xi_n\|+|c|^2 \to |c|^2 .
\]
However, this contradicts $|c|<1$. Therefore, $E$ must be finite-rank. Now, let $W_{E(\cH)}(C) = \{ \langle C\xi,\xi \rangle : \xi \in E(\cH),~\|\xi\|=1\}$ for each $C \in \cB (\cH)$. Since $W_{E(\cH)}(C) = W_{E(\cH)}(ECE)$ and $ECE$ can be viewed as an operator on the finite-dimensional subspace $E(\cH)$ of $\cH$, it turns out that $W_{E(\cH)}(C)$ is always closed. Now, recall that $\rho ((B+dI)^*(A+cI))=0$, which implies that $0 \in W_{E(\cH)}((B+dI)^*(A+cI))$ by the Toeplitz-Hausdorff theorem. This proves $A+cI \perpq B+dI$.

Suppose that $|c|\geq 1$. It follows from
\[
1 \leq |c|^2 =\rho_\infty (|A+cI|^2) \leq 1
\]
that $\rho_\infty (|A+cI|)=|c|=1=\|A+cI\|$. If $\rho_\infty ((B+dI)^*(A+cI)) = c\overline{d} = 0$, then $A+cI \perpq B+dI$. We assume that $c\overline{d} \neq 0$. In this case, $\rho ((B+dI)^*(A+cI)) =0$ implies $\rho_u ((B+dI)^*(A+cI)) =-c\overline{d}\|\rho_s\|$. Meanwhile, if $C_1,\ldots ,C_n \in \cK (\cH)$, then $C_1^*,\ldots, C_n^* \in \cK (\cH)$ has the separable range and so is $\bigvee_{j=1}^n R(C_j^*)$ by $(\bigvee_{j=1}^n R(C_j^*))(\cH) = [\bigcup_{j=1}^n \overline{C_j^*(\cH)}]$. From this, we derive $\bigwedge_{j=1}^n N(C_j) =\bigwedge_{j=1}^n (I-R(C_j^*)) = I-\bigvee_{j=1}^nR(C_j^*) \neq 0$, where $R(C_j^*)$ and $N(C_j)$ are respectively the projections onto $\overline{C_j^*(\cH)}$ and $\ker C_j$. Now, let $F = N(A) \wedge N(B) \neq 0$, and let $\zeta$ be a unit vector in $F(\cH)$. Set $E' = E \vee E_\zeta$ for short. Then, we get
\[
0 =\langle (B^*A+cB+\overline{d}A)\zeta, \zeta \rangle \in W_{E'(\cH)} (B^*A+cB+\overline{d}A)
\]
where $W_{E'(\cH)}(C) = \{ \langle C\xi,\xi \rangle :\xi \in E'(\cH),~\|\xi\|=1\}$. Since $E'(B^*A+cB+\overline{d}A)E' \in \cK (E'(\cH))$, \cite[Theorem 1]{BGS72} ensures that $W_{E'(\cH)}(B^*A+cB+\overline{d}A)$ is closed, which in turn implies that
\begin{align*}
W_{E'(\cH)}((B+dI)^*(A+cI)) 
&= W_{E'(\cH)}(B^*A+cB+\overline{d}A+c\overline{d}I) \\
&= \{ \lambda +c\overline{d} : \lambda \in W_{E'(\cH)}(B^*A+cB+\overline{d}A)\}
\end{align*}
is also closed. Hence,
\[
-c\overline{d}\|\rho_s\|\|\rho_u\|^{-1} = (\|\rho_u\|^{-1}\rho_u )((B+dI)^*(A+cI)) \in W_{E'(\cH)}((B+dI)^*(A+cI)).
\]
From this and $c\overline{d}=\omega_{\zeta} ((B+dI)^*(A+cI)) \in W_{E'(\cH)}((B+dI)^*(A+cI))$, the Toeplitz-Hausdorff theorem generates a unit vector $\eta \in E'(\cH)$ such that $\omega_\eta ((B+dI)^*(A+cI))=0$. Finally, we remark that $|A+cI|\zeta =\zeta$ by
\[
\langle |A+cI|\zeta ,\zeta \rangle \geq \langle |A+cI|^2\zeta ,\zeta \rangle = \langle (A^*A+cA^*+\overline{c}A+|c|^2 )\zeta ,\zeta \rangle =|c|^2=1 ,
\]
which ensures $|A+cI|\eta =\eta$. Thus, $A+cI \perpq B+dI$.
\end{proof}
\end{theorem}
\begin{corollary}\label{Cor:cpt-op}

Let $\cH$ be an infinite-dimensional Hilbert space. Then, $\perpq$ and $\perpn$ are equivalent on $\cK (\cH)^\sim$ if and only if $\cH$ is nonseparable.
\end{corollary}
Recall that a $C^*$-algebra $\A$ is said to be \emph{elementary} if $\A$ is $*$-isomorphic to some $\cK (\cH)$, \emph{simple} if $\A$ contains no nontrivial closed two-sided ideal, and \emph{primitive} if there exists a faithful irreducible representation of $\A$. It is obivious that an elementary $C^*$-algebra is simple. Since the kernel of each representation is a closed two-sided ideal, a simple $C^*$-algebra is primitive. It is known that the primitivity of a $C^*$-algebra implies primeness; see, for example,~\cite[II.5.4.6]{Bla06}. By Corollary~\ref{Cor:cpt-op}, we have a non-simple $C^*$-algebra in which $\perpq$ and $\perpn$ are equivalent, and a primitive $C^*$-algebra in which $\perpq$ and $\perpn$ are not equivalent. This situation for $C^*$-algebras is summarized as follows:
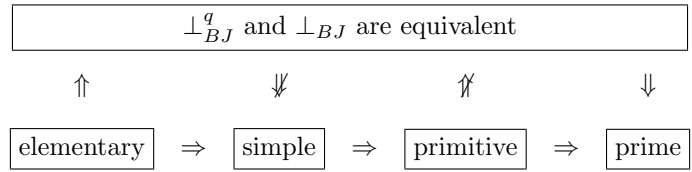
\begin{align*}
\begin{array}{c}
\text{\fbox{\hspace{6.125em}$\perpq$ and $\perpn$ are equivalent\hspace{6.125em}}}\\[1em]
\begin{array}{ccccccc}
\Uparrow & & \not \Downarrow & & \not \Uparrow & & \Downarrow\\[1em]
\text{\fbox{elementary}} & \Rightarrow & \text{\fbox{simple}} & \Rightarrow & \text{\fbox{primitive}} & \Rightarrow & \text{\fbox{prime}}
\end{array}
\end{array}
\end{align*}

The following problems are still open.
\begin{problem}

Let $\A$ be a simple $C^*$-algebra. Are $\perpq$ and $\perpn$ equivalent on $\A$?
\end{problem}

\begin{problem}

Let $\A$ be a $C^*$-algebra on which $\perpq$ and $\perpn$ are equivalent. Is $\A$ primitive?
\end{problem}



\end{document}